\documentclass{revtex4}%
\usepackage{amsfonts,pstricks}
\usepackage{amsmath}
\usepackage{amssymb}
\usepackage{graphicx}
\usepackage{amsfonts}%
\newtheorem{theorem}{Theorem}

\newtheorem{observation}[theorem]{Observation}

\newtheorem{lemma}[theorem]{Lemma}

\newenvironment{proof}[1][Proof]{\noindent\textbf{#1.} }{\ \rule{0.5em}{0.5em}}
\begin{document}
\title[Counting paths in Bratteli diagrams for $SU(2)_{k}$]{Counting paths in Bratteli diagrams for $SU(2)_{k}$}
\author{Toufik Mansour}
\affiliation{Department of Mathematics, University of Haifa, Haifa 31905, Israel}
\author{Simone Severini}
\affiliation{Institute for Quantum Computing and Department of Combinatorics \&
Optimization, University of Waterloo, Waterloo N2L 3G1, Canada}

\begin{abstract}
It is known that the Hilbert space dimensionality for quasiparticles in an
$SU(2)_{k}$ Chern-Simons-Witten theory is given by the number of directed
paths in certain Bratteli diagrams. We present an explicit formula for these
numbers for arbitrary $k$. This is on the basis of a relation with Dyck paths
and Chebyshev polynomials.

\end{abstract}
\maketitle

\section{Introduction}

The so-called Bratteli diagrams have been introduced by Bratteli in 1972
\cite{br} for the classification of some classes of $C^{\ast}$-algebras. For
our scope, it is sufficient to give here an \textquotedblleft algebra
free\textquotedblright\ description of this notion (see, \emph{e.g.},
\cite{be} or \cite{d}). More specifically, lattice points will be labeled by
numbers instead of Young diagrams. A \emph{Bratteli diagram} for $SU(2)_{k}$
is defined as a finite digraph $D_{k}=(V,E)$, where:
\begin{itemize}
\item The vertices of $D_{k}$ are associated to lattice points in the positive
quadrant of the Cartesian plane:
$$(j,i) \in(\mathbb{Z}^{\geq0})^{\times2}\mbox{ with }j\geq i.$$
The vertex set $V$ is $\{(j,i) :0\leq i\leq k,j$; $i$ and $j$ have
the same parity$\}$. Here $0$ is assumed to be even.

\item There is an arc $((j,i),(j^{\prime},i^{\prime}))$ in $E$ if $j^{\prime}=j+1$
and $i^{\prime}=i\pm1$.

\item A \emph{directed path} from a vertex $(j,i)$ to a vertex $(j^{\prime
},i^{\prime})$ of length $n$ is a sequence of $n$ arcs of the form
$(j,i)(j+1,l_{1}),(j+1,l_{1})(j+2,l_{2}),\ldots,(j+n-1,l_{n-1})(j^{\prime
},i^{\prime})$. Let $D_{k}(x,y)$ be the number of directed paths from the
vertex $(0,0)$ to the vertex $(j,i)$ in the Bratteli diagram $D_{k}$. In a
Bratteli diagram $D_{k}$, the vertex $(j,i)$ is labeled by the number
$D_{k}(i,j)$. In our notation, the function $f$ introduced in the general
definition of a Bratteli diagram is exactly $D_{k}(i,j)$. Notice that
$$D_{k}(i,j)=D_{k}(i-1,j-1)+D_{k}(i+1,j-1)$$
with the initial conditions
$D_{k}(0,j)=D_{k}(1,j-1)$ and $D_{k}(k,j)=D_{k}(k-1,j-1)$.
\end{itemize}

Let $$d^{+}(i)=\left\vert \left\{  j:(i,j)\in A(D_{k})\right\}  \right\vert $$
and $$d^{-}(i)=\left\vert \left\{  j:(j,i)\in A(D_{k})\right\}  \right\vert $$
be the \emph{indegree} and the \emph{outdegree} of a vertex $i$, respectively.
Notice that $d^{+}(0,0)=1$, $d^{-}(0,0)=0$, and $d^{\pm}(0,j)=d^{\pm
}(k,i)=d^{\pm}(l,l)=1$. All other vertices have indegree and outdegree $2$.
The figure below illustrates the Bratteli diagram $D_{2}$.

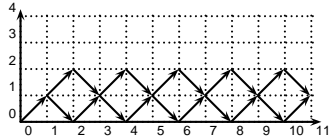
\begin{figure}[h]
\begin{pspicture}(0,0)(3.85,2.5)
\psline{->}(0,0)(4,0) \psline{->}(0,0)(0,1.5)
\psgrid[xunit=10pt,yunit=10pt,subgriddiv=1,griddots=5,gridlabels=5pt](0,0)(11,4)
\multips(0,0)(.7,0){5}{\psline{->}(0,0)(.35,.35)\psline{->}(.35,.35)(.7,0)}
\multips(.35,.35)(.7,0){5}{\psline{->}(0,0)(.35,.35)\psline{->}(.35,.35)(.7,0)}
\psline{->}(3.5,0)(3.85,.35)
\end{pspicture}\caption{Bratteli diagram $D_{2}$.}%
\end{figure}

The number $D_{k}(i,j)$ corresponds to the Hilbert space dimensionality for
$i$ $q$-spin $1/2$ quasiparticles having total $q$-spin $j$, in an $SU(2)_{k}$
Chern-Simons-Witten theory (see, \emph{e.g.}, \cite{ho} and the references
contained therein). An asymptotic expression and a recurrence relation for
$D_{k}(i,j)$ have been pointed out by Slingerland and Bais (see \cite[Section
2.4]{sb}). More precisely, they showed that
\begin{equation}
D_{k}(i,j)\approx\left(  2\cos\left(  \frac{\pi}{k+2}\right)  \right)  ^{j},
\label{eqs0}%
\end{equation}
where $i+j=0(\mod2)$ (From the definitions we have that $D_{k}(i,j)=0$ for all
$i+j=1(\operatorname{mod}2).$) We present here an explicit formula for
$D_{k}(i,j)$, which implies the asymptotic expression in \eqref{eqs0}. In
particular, we find an explicit formula for the total number of directed paths
from the vertex $\left(  0,0\right)  $ to the vertex $(i,j)$ in the Bratteli
diagram $D_{k}$, for any $k$, as shown in the next section. Our tools will be
Dyck paths and Chebyshev polynomials.

\section{Main result}

\emph{Chebyshev polynomials of the second kind}\textit{ }are defined by
\[
U_{r}(\cos\theta)=\frac{\sin\left(  (r+1)\theta\right)  }{\sin\theta},
\]
for $r\geq0$. Evidently, $U_{r}(x)$ is a polynomial of degree $r$ in $x$ with
integer coefficients. Chebyshev polynomials were invented for the needs of
approximation theory, but are also widely used in various other branches of
mathematics, including algebra, number theory, and the study of lattice paths
in combinatorics (see \cite{Ri}). For $k\geq0$, we define $R_{k}(x)$ by
\[
R_{k}(x)=\frac{U_{k-1}\left(  \frac{1}{2\sqrt{x}}\right)  }{\sqrt{x}%
U_{k}\left(  \frac{1}{2\sqrt{x}}\right)  }.
\]
For example, $R_{0}(x)=0$, $R_{1}(x)=1$, and $R_{2}(x)=1/(1-x)$. It is easy to
see that for any $k$, $R_{k}(x)$ is a rational function in $x$.

A \emph{Dyck path} is a lattice path in the plane integer lattice
$(\mathbb{Z}^{\geq0})^{\times2}$ consisting of \emph{up-steps} $u=(1,1)$ and
\emph{down-steps} $d=(1,-1)$. It follows that a Dyck path never passes below
the $x$-axis. The length of a Dyck path is defined as the number of its
up-steps and down-steps. From the definitions we can state the following observation.

\begin{observation}
The number $D_{k}(i,j)$ is exactly the number of Dyck paths below the line
$y=k+1$, starting at the origin and ending at the point with $x$-coordinate
$j$ and $y$-coordinate $i$.
\end{observation}

Each Dyck path $P$ below the line $y=k+1$, starting at the origin and ending
at $(j,i)$, has the following form:
\[
P=P_{1}uP_{2}uP_{3}\cdots uP_{i+1},
\]
where $P_{s}$ is a Dyck path of height at most $k+1-s$, and the length of $P$
is exactly $j$ (each up-step and down-step is counted as a unit step).

Let us fix a variable $x$ for the generating function for counting the number
of up-steps and down-steps in a Dyck path. Using the fact that the generating
function for the number of Dyck paths of length $2n$ with height at most $s$
is given by $R_{s+1}(x^{2})$ (see \cite{Ma}), we get that the generating
function
\[
D_{k}(x;i)=\sum_{j\geq0}D_{k}(i,j)x^{j}%
\]
is given by
\begin{equation}
D_{k}(x;i)=R_{k+1}(x^{2})\prod_{r=1}^{i}(xR_{k+1-r}(x^{2})), \label{eqss}%
\end{equation}
which is equivalent to
\begin{equation}
\label{eqssaa}D_{k}(x;i)=x^{i}\prod_{r=0}^{i}\frac{U_{k-r}\left(  \frac{1}%
{2x}\right)  }{xU_{k+1-r}\left(  \frac{1}{2x}\right)  }=\frac{U_{k-i}\left(
\frac{1}{2x}\right)  }{xU_{k+1}\left(  \frac{1}{2x}\right)  }.
\end{equation}
Using the fact that the roots of $U_{k}(x)$ are $$\cos\left(  \frac{r\pi}{k+1}\right),\ r=1,2,\ldots,k$$
we obtain that the minimal positive pole of the function $D_{k}(x;i)$ is given
by $\left(  2\cos\left(  \frac{\pi}{k+2}\right)  \right)  ^{-1}$. Therefore,
the asymptotic behaviour of the function $D_{k}(i,j)$ is given by
\[
D_{k}(i,j)\approx\left(  2\cos\left(  \frac{\pi}{k+2}\right)  \right)  ^{j},
\]
where $i+j=0(\operatorname{mod}2)$, as described in \cite[Equation 26]{sb}.

In order to find an explicit formula for $D_{k}(i,j)$ we need the following lemma.

\begin{lemma}
\label{lema0} The generating function $$\frac{U_{k-i}(x)}{U_{k+1}(x)}$$ is given
by
\[
\frac{1}{k+2}\sum_{r=1}^{k+1}\frac{(-1)^{r+1}U_{k-i}(\rho_{k+1,r})\sin
^{2}\frac{\pi r}{k+2}}{x-\rho_{k+1,r}},
\]
where $\rho_{m,r}=\cos\bigl(\frac{r\pi}{m+1}\bigr)$.
\end{lemma}

\begin{proof}
Let us compute the partial fraction decomposition of $$\frac{U_{k-i}%
(x)}{U_{k+1}(x)}.$$ By general principles, it is $$\sum_{r=1}^{k+1}%
\frac{a_{k+1,r}}{x-\rho_{k+1,r}},$$ where $\rho_{m,r}=\cos\bigl(\frac{r\pi
}{m+1}\bigr)$ are the zeros of the $m$-th Chebyshev polynomials of the second
kind. Now, $$a_{k+1,r}=\frac{U_{k-i}(\rho_{k+1,r})}{\left.  \frac{d}%
{dx}U^{\prime}_{k+1}(x)\right|  _{\rho_{k+1,r}}}$$
and note that
\[
U^{\prime}_{k+1}(x)=\frac{dU_{k+1}(x)}{dx}=\frac d{d\theta}\Big(\frac
{\sin(k+2)\theta}{\sin\theta}\Big)\cdot\frac{d\theta}{dx}.
\]
We work out that
\[
\frac{dU_{k+1}}{d\theta}=\frac{(k+2)\cos(k+2)\theta\cdot\sin\theta
-\sin(k+2)\theta\cdot\cos\theta}{\sin^{2}\theta},
\]
and if we plug in $x=\rho_{k+1,r}$ simplification occurs, since certain terms
are just zero; we obtain that
\begin{align*}
&  \frac{dU_{k+1}}{d\theta}(\arccos\rho_{k+1,r})\\
&  =\frac{(k+2)\cos(k+2)\theta\cdot\sin\theta- \sin(k+2)\theta\cdot\cos\theta
}{\sin^{2}\theta}\bigg|_{\theta=\frac{\pi r}{k+2}}\\
&  =\frac{(k+2)\cos(\pi r)\cdot\sin\frac{\pi r}{k+2}- \sin(\pi r)\cdot
\cos\frac{\pi r}{k+2}}{\sin^{2}\frac{\pi r}{k+2}}\\
&  =\frac{(k+2)\cos(\pi r)\cdot\sin\frac{\pi r}{k+2}}{\sin^{2}\frac{\pi
r}{k+2}}\\
&  =\frac{(k+2)(-1)^{r}}{\sin\frac{\pi r}{k+2}}.
\end{align*}
Further $\frac{dx}{d\theta}=-\sin\theta$, so together
\begin{align*}
a_{k+1,r}  &  =\frac{U_{k-i}(\rho_{k+1,r})}{U^{\prime}_{k+1}(\rho_{k+1,r})}\\
&  =\frac{1}{k+2}(-1)^{r+1}U_{k-i}(\rho_{k+1,r})\sin^{2}\frac{\pi r}{k+2},
\end{align*}
which completes the proof.
\end{proof}

Now we are ready to give an explicit formula for $D_{k}(i,j)$.

\begin{theorem}
For all $k,i,j$, the number $D_{k}(i,j)$ of directed paths from the vertex
$\left(  0,0\right)  $ to the vertex with $x$-coordinate $j$ and
$y$-coordinate $i$ in the Bratteli diagram $D_{k}$ is given by
\[
\frac{2}{k+2}\sum_{r=1}^{k+1}(-1)^{r+1}U_{k-i}(\rho_{k+1,r})\sin^{2}\left(
\frac{\pi r}{k+2}\right)  (2\rho_{k+1,r})^{j},
\]
where $\rho_{m,r}=\cos\bigl(\frac{r\pi}{m+1}\bigr)$.
\end{theorem}

\begin{proof}
Lemma~\ref{lema0} and Equation \ref{eqssaa} give that the generating function
$D_{k}(x;i)$ is given by
\[
\frac{2}{k+2}\sum_{r=1}^{k+1}\frac{(-1)^{r+1}U_{k-i}(\rho_{k+1,r})\sin
^{2}\left(  \frac{\pi r}{k+2}\right)  }{1-2\rho_{k+1,r}x},
\]
which implies that the $x^{j}$ coefficient in the generating function
$D_{k}(x;i)$ is given by
\[
\frac{2}{k+2}\sum_{r=1}^{k+1}(-1)^{r+1}U_{k-i}(\rho_{k+1,r})\sin^{2}\left(
\frac{\pi r}{k+2}\right)  (2\rho_{k+1,r})^{j},
\]
as claimed.
\end{proof}

For example, the above theorem for $k=1,2,3,4,5$ gives the following values:

\begin{itemize}
\item $D_{1}(0,2j)=1$ and $D_{1}(1,2j+1)=1$ (otherwise $D_{1}(i,j)=0$).
\begin{figure}[h]
\begin{pspicture}(0,0)(3.85,1.6)
\psline{->}(0,0)(4,0) \psline{->}(0,0)(0,1.5)
\psgrid[xunit=10pt,yunit=10pt,subgriddiv=1,griddots=5,gridlabels=0pt](0,0)(11,4)
\multips(0,0)(.7,0){5}{\psline[linewidth=.01]{->}(0,0)(.35,.35)\psline[linewidth=.01]{->}(.35,.35)(.7,0)}
\multips(.35,.35)(.7,0){5}{\psline[linewidth=.01]{->}(0,0)(.35,.35)\psline[linewidth=.01]{->}(.35,.35)(.7,0)}
\psline[linewidth=.01]{->}(3.5,0)(3.85,.35) \put(-.2,0){\tiny
1}\put(.1,.4){\tiny 1}\put(.6,.75){\tiny 1}\put(.6,-.28){\tiny
1}\put(1,.5){\tiny 2} \put(1.35,.75){\tiny 2}\put(1.3,-.28){\tiny
2}\put(1.7,.5){\tiny 4} \put(2,.75){\tiny 4}\put(2,-.28){\tiny
4}\put(2.4,.5){\tiny 8} \put(2.7,.75){\tiny 8}\put(2.7,-.28){\tiny
8}\put(3.05,.5){\tiny 16} \put(3.4,.75){\tiny
16}\put(3.4,-.28){\tiny 16}\put(3.7,.5){\tiny 32}
\end{pspicture}\caption{Number of paths in the Bratteli diagram $D_{2}$.}%
\label{fig2}%
\end{figure}
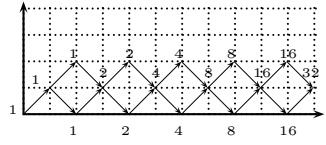

\item $D_{2}(0,2j)=D_{2}(2,2j)=2^{j-1}$ and $D_{2}(1,2j+1)=2^{j}$, see
Figure~\ref{fig2}.\newline

\item $D_{3}(0,2j)=F_{2j-1}$, $D_{3}(1,2j)=F_{2j+1}$, $D_{3}(2,2j+2)=F_{2j+2}$
and $D_{3}(3,2j+3)=F_{2j+2}$, where $F_{m}$ is the $m$-th Fibonacci number
(defined as $F_{0}=0$, $F_{1}=1$, and $F_{n}=F_{n-1}+F_{n-2}$ for any $n$),
see Figure~\ref{fig3}. It was pointed out in \cite{ho} that the theory of
$q$-spin $1$ quasiparticles in $SU(2)_{3}$ is equivalent to $SO(3)_{3}$, a
theory which is in fact also known as the Fibonacci anyon theory \cite{pr}.
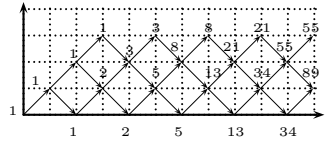
\begin{figure}[h]
\begin{pspicture}(0,0)(3.85,2)
\psline{->}(0,0)(4,0) \psline{->}(0,0)(0,1.5)
\psgrid[xunit=10pt,yunit=10pt,subgriddiv=1,griddots=5,gridlabels=0pt](0,0)(11,4)
\multips(0,0)(.7,0){5}{\psline[linewidth=.01]{->}(0,0)(.35,.35)\psline[linewidth=.01]{->}(.35,.35)(.7,0)}
\multips(.7,.7)(.7,0){4}{\psline[linewidth=.01]{->}(0,0)(.35,.35)\psline[linewidth=.01]{->}(.35,.35)(.7,0)}
\multips(.35,.35)(.7,0){5}{\psline[linewidth=.01]{->}(0,0)(.35,.35)\psline[linewidth=.01]{->}(.35,.35)(.7,0)}
\psline[linewidth=.01]{->}(3.5,0)(3.85,.35)\psline[linewidth=.01]{->}(3.5,.7)(3.85,1.05)
\put(-.2,0){\tiny 1}\put(.1,.4){\tiny 1}\put(.6,.75){\tiny
1}\put(.6,-.28){\tiny 1}\put(1,.5){\tiny 2} \put(1,1.1){\tiny 1}
\put(1.35,.8){\tiny 3}\put(1.3,-.28){\tiny 2} \put(1.7,.5){\tiny
5}\put(1.7,1.1){\tiny 3} \put(1.95,.85){\tiny 8}\put(2,-.28){\tiny
5} \put(2.4,.5){\tiny 13}\put(2.4,1.1){\tiny 8}
\put(2.65,.85){\tiny 21}\put(2.7,-.28){\tiny 13}
\put(3.05,.5){\tiny 34}\put(3.05,1.1){\tiny 21}
\put(3.35,.85){\tiny 55}\put(3.4,-.28){\tiny 34}
\put(3.7,.5){\tiny 89}\put(3.7,1.1){\tiny 55}
\end{pspicture}\caption{Number of paths in the Bratteli diagram $D_{3}$.}%
\label{fig3}%
\end{figure}

\item $D_{4}(0,2j+2)=D_{4}(1,2j+1)=(3^{j}+1)/2$,\newline$D_{4}(2,2j+2)=3^{j}$,
and\newline$D_{4}(3,2j+3)=D_{4}(4,2j+4)=(3^{j+1}-1)/2$,\newline for all
$j\geq0$.\newline

\item $D_{5}(0,2j)=a_{j}$, $D_{5}(1,2j+1)=a_{j+1}$,\newline$D_{5}%
(2,2j+2)=a_{j+2}-a_{j+1}$,\newline$D_{5}(3,2j+3)=a_{j+3}-2a_{j+2}$%
,\newline$D_{5}(4,2j+4)=D_{5}(5,2j+5)=a_{j+4}-3a_{j+3}+a_{j+2}$,\newline for
all $j\geq0$, where $a_{j}$ is a sequence satisfies the following recurrence
relation $a_{m}=5a_{m-1}-6a_{m-2}+a_{m-3}$ with the initial conditions
$a_{0}=a_{1}=1$ and $a_{2}=2$. Note that, the generating function for the
sequence $a_{n}$ is given by $$\frac{1-4x+3x^{2}}{1-5x+6x^{2}-x^{3}}.$$
\end{itemize}

Using the continued fraction representation of the generating function
$R_{k}(x)$, namely
\[
R_{k}(x)=\frac{1}{1-\dfrac{x}{1-\frac{x}{\ddots}}}%
\]
with exactly $k$-levels (see~\cite[Lemma 3.1]{MV}), we obtain that
\[
R_{\infty}(x)=\lim_{k\rightarrow\infty}R_{k}(x)=C(x):=\frac{1-\sqrt{1-4x}}%
{2x}.
\]
Hence, From \eqref{eqss} we obtain that the generating function
\[
D_{\infty}(x;i)=\lim_{k\rightarrow\infty}D_{k}(x;i)
\]
is given by
\[
D_{\infty}(x;i)=x^{i}C^{i+1}(x^{2}).
\]
The generating function $C(x)$ satisfies $$C(x)=1+xC^{2}(x),$$ thus by the
Lagrange inversion formula \cite[Sec. 5.4]{stanley} we obtain that
\[
D_{\infty}(i,j)=\frac{i+1}{j+1}\binom{j+1}{\frac{j-i}{2}},
\]
where $i+j=0(\operatorname{mod}2)$, as shown in \cite[Equation 21]{sb}.

\end{document}